\documentclass[11pt]{article}

\usepackage{amsmath,amssymb,latexsym}
\usepackage{MnSymbol}
\usepackage[latin1]{inputenc}
\usepackage[dvips]{graphicx}
\usepackage{hyperref}
\usepackage{color}
\usepackage{mathrsfs}
\usepackage{enumerate}
\usepackage{tikz}
\usepackage{xifthen}
\usepackage{verbatim}
\usepackage{soul}

\newtheorem{theorem}{Theorem}[section]

\newtheorem{lemma}[theorem]{Lemma}
\newtheorem{proposition}[theorem]{Proposition}

\newtheorem{corollary}[theorem]{Corollary}
\newtheorem{remark}[theorem]{Remark}

\newtheorem{problem}[theorem]{Problem}

\newcommand{\proof}{\noindent{\bf Proof.\ }}
\newcommand{\qed}{\hfill $\square$ \bigskip}

\newcommand{\gp}{{\rm gp}}
\newcommand{\sgp}{{\rm sgp}}
\newcommand{\sjc}{{\rm sjc}}

\begin{document}

\title{A Steiner general position problem in graph theory}

\author{
Sandi Klav\v zar $^{a,b,c}$
\and
Dorota Kuziak $^{d}$
\and
Iztok Peterin $^{e,c}$
\and
Ismael G. Yero $^d$
}

\date{}

\maketitle

\begin{center}
$^a$ Faculty of Mathematics and Physics, University of Ljubljana, Slovenia \\
\medskip

$^b$ Faculty of Natural Sciences and Mathematics, University of Maribor, Slovenia \\
\medskip

$^{c}$ Institute of Mathematics, Physics and Mechanics, Ljubljana, Slovenia \\
{\tt sandi.klavzar@fmf.uni-lj.si}
\medskip

$^d$ Departamento de Matem\'aticas, Universidad de C\'adiz, Algeciras, Spain \\
{\tt dorota.kuziak@uca.es, ismael.gonzalez@uca.es}
\medskip

$^{e}$ Faculty of Electrical Engineering and Computer Science, University of Maribor, Slovenia \\
{\tt iztok.peterin@um.si}
\end{center}

\maketitle

\begin{abstract}
Let $G$ be a graph. The Steiner distance of $W\subseteq V(G)$ is the minimum size of a connected subgraph of $G$ containing $W$. Such a subgraph is necessarily a tree called a Steiner $W$-tree. The set $A\subseteq V(G)$ is a $k$-Steiner general position set if $V(T_B)\cap A = B$ holds for every set $B\subseteq A$ of cardinality $k$, and for every Steiner $B$-tree $T_B$. The $k$-Steiner general position number $\sgp_k(G)$ of $G$ is the cardinality of a largest $k$-Steiner general position set in $G$. Steiner cliques are introduced and used to bound $\sgp_k(G)$ from below. The $k$-Steiner general position number is determined for trees, cycles and joins of graphs. Lower bounds are presented for split graphs, infinite grids and  lexicographic products. The lower bound for the latter products leads to an exact formula for the general position number of an arbitrary lexicographic product.  
\end{abstract}


\medskip\noindent
{\bf AMS Subj.\ Class.\ (2010)}: 05C12, 05C05, 05C76

%%%%%%%%%%%%%%%%%%%%%%%%%%%%%%%%%%%%%%%%%%%%%%%%%%%%%
%%%%%%%%%%%%%%%%%%%%%%%%%%%%%%%%%%%%%%%%%%%%%%%%%%%%%
\section{Introduction}
\label{sec:intro}
%%%%%%%%%%%%%%%%%%%%%%%%%%%%%%%%%%%%%%%%%%%%%%%%%%%%%
%%%%%%%%%%%%%%%%%%%%%%%%%%%%%%%%%%%%%%%%%%%%%%%%%%%%%

In this work, $G = (V(G), E(G))$ denotes a simple graph. The \emph{distance} $d_G(u,v)$ between two vertices $u$ and $v$ of $G$ is the minimum number of edges on a $u,v$-path in $G$. If there is no such path, then we set $d_G(u,v)=\infty$. A $u,v$-path of length $d_G(u,v)$ is called a {\em $u,v$-geodesic}.

A {\em general position set} of a graph $G$ is a set of vertices $S\subseteq V(G)$ such that no three vertices from $S$ lie on a common geodesic. The cardinality of a largest possible general position set is the {\em general position number} $\gp(G)$ of $G$. The problem of finding the general position number was independently introduced in~\cite{manuel-2018a, ullas-2016} and earlier studied on hypercubes in~\cite{korner-1995}. The paper~\cite{manuel-2018a} opened a wide interest for the topic, articles~\cite{anand-2019, ghorbani-2019, klavzar-2019+, klavzar-2021, klavzar-2019, manuel-2018b, neethu-2021, patkos-2020, thomas-2020, tian-2020, tian-2021} bring many different results on the general position sets and numbers. General position sets were also generalized to general $d$-position sets, where $d$ is a threshold on the length of geodesics on which triples of vertices are not allowed to lie~\cite{klavzar-2021+}.

For a nonempty set $W\subseteq V(G)$, the {\em Steiner distance} of $W$, denoted by $d_G(W)$, is the minimum size of a connected subgraph of $G$ containing $W$~\cite{Chartrand-1989}. Such a subgraph is clearly a tree called a \emph{Steiner $W$-tree}. If $G$ is not connected and the vertices of $W$ lie in at least two components of $G$, then no Steiner $W$-tree exists and we set $d_G(W) = \infty$. Papers~\cite{dankelmann-2021, li-2016, martinez-2018, nielsen-2009, zhang-2019} represent a selection of  studies on Steiner trees.

We now introduce the key new concept. Let $k\in {\mathbb N}$ and let $G$ be a graph. Then $A\subseteq V(G)$ is a {\em $k$-Steiner general position set} if for every set $B\subseteq A$ of cardinality $k$ (from now on a $k$-set), and for every Steiner $B$-tree $T_B$, it follows that $V(T_B)\cap A = B$. In other words, $A$ is a $k$-Steiner general position set if no $k+1$ distinct vertices from $A$ lie on a common Steiner $B$-tree, where $B\subseteq A$ and $|B| = k$. Clearly, if $|A|\le k$, then $A$ is $k$-Steiner general position set. Hence we may define the {\em $k$-Steiner general position number} of $G$, denoted by $\sgp_k(G)$, as the cardinality of a largest $k$-Steiner general position set in $G$. A $k$-Steiner general position set of cardinality $\sgp_k(G)$ will be called a {\em $k$-sgp-set}. Note that $\sgp_2(G) = \gp(G)$.

Recall that if $G$ is a graph and $A\subseteq V(G)$, then $A$ is {\em Steiner convex} if for any subset $B\subseteq A$, all vertices in every Steiner tree $T_B$ belong to $A$, see~\cite{caceras-2008, gologranc-2015, gologranc-2018}. The new concept of the Steiner general position set is hence a concept dual to the Steiner convex set.

We proceed as follows. In the rest of this section some further definitions are listed. In the next section we present basic bounds for $\sgp_k(G)$ and settle the extreme case when $\sgp_k(G)$ is equal to the order of $G$. The section devoted to $\sgp_k(G)$ for trees and cycles follows. In the fourth section we present an exact result for $\sgp_k(G\vee H)$, where $G\vee H$ represents the join of graphs $G$ and $H$. This enables us to present several exact results for some known families of graphs. After that comes a section with the $k$-Steiner general position number of lexicographic products. The before last section brings lower bounds for $\sgp_k(G)$ for split graphs and infinite grid graphs. We present several open problems and questions in the last section.

We conclude the introduction by some necessary terminology and notation. By $n(G)$ we denote the order of a graph $G$. As usual, $\delta(v)$ represents the \emph{degree} of a vertex $v\in V(G)$, \emph{i.e.}, the number of neighbors of $v$. If $T$ is a tree, then $L(T)$ is the set of its leaves and $\ell(T)=|L(T)|$. The \emph{clique number} $\omega(G)$ of $G$ is the cardinality of a largest complete subgraph of $G$. If $A\subseteq V(G)$, then the subgraph of $G$ induced by $A$ is denoted by $G[A]$. A graph $G$ is $d$-\emph{connected}, where $1\le d < n(G)$, if the removal of fewer than $d$ vertices from $G$ always yields a connected graph. By $\overline{G}$ we denote the \textit{complement graph} of $G$; it is defined by $V(\overline{G})=V(G)$ and $uv\in E(\overline{G})$ if and only if $uv\notin E(G)$. For positive integers $i < j$ we use the notation $[i:j] = \{i,i+1,\ldots ,j\}$. Other definitions that will be needed, like the join of graphs, the lexicographic product and others, will be given along the way. 

%%%%%%%%%%%%%%%%%%%%%%%%%%%%%%%%%%%%%%%%%%%%%%%%%%%%%
%%%%%%%%%%%%%%%%%%%%%%%%%%%%%%%%%%%%%%%%%%%%%%%%%%%%%
\section{Bounding the $k$-Steiner general position number}
\label{sec:bounding}
%%%%%%%%%%%%%%%%%%%%%%%%%%%%%%%%%%%%%%%%%%%%%%%%%%%%%
%%%%%%%%%%%%%%%%%%%%%%%%%%%%%%%%%%%%%%%%%%%%%%%%%%%%%

To bound the $k$-Steiner general position number from below, we introduce in this section $k$-Steiner cliques, a concept that might be of interest also elsewhere. We characterize graphs $G$ with $\sgp_k(G) = n(G)$ and discuss monotonicity of $\sgp_k(G)$ with respect to the parameter $k$.

Let $G$ be a connected graph. Since every nonempty set $A\subseteq V(G)$ is trivially a $1$-Steiner general position set and because $\sgp_{n(G)}(G) = n(G)$, in the rest of the paper we restrict our considerations to the $k$-Steiner general position sets of $G$ with $k\in [2:n(G)-1]$.

Induced subgraphs of cliques are cliques, hence $\sgp_k(G) \ge \omega(G)$. If $k\ge 2$ is a fixed integer, then $A\subseteq V(G)$ is a {\em $k$-Steiner clique} if $G[B]$ is connected for every $k$-set $B\subseteq A$. The cardinality of a largest $k$-Steiner clique will be denoted by $s\omega_k(G)$. Note that $s\omega_2(G)=\omega(G)$. We make the assumption that every set on $k-1$ or less vertices is a $k$-Steiner clique, since indeed there are no $k$-sets as subsets of such set, and we can use this fact to deal with disconnected graphs as well. If $G$ is not connected, and the order of each component of $G$ is smaller than $k$, then every set on $\min\{n(G),k-1\}$ vertices of $G$ represents a $k$-Steiner clique. Hence, in such a case we have $s\omega_k(G)=\min \{n(G),k-1\}$.

Notice that if $G$ is a connected graph with $n(G)\ge k$, then $s\omega_k(G)\geq k$. Further, a clique is a $k$-Steiner clique for every $k\ge 2$, thus $s\omega_k(G) \ge \omega(G)$. Moreover, if $A$ is a $k$-Steiner clique of $G$, then it is also a $k$-Steiner general position set, hence $\sgp_k(G) \ge  s\omega_k(G)$. This discussion can be summarized as follows.

\begin{remark}
\label{rem:trivial}
If $G$ is a connected graph and $k \in [2:n(G)-1]$, then
$$\max\{k,\omega(G)\}\le s\omega_k(G)\le \sgp_k(G)\le n(G)\,.$$
\end{remark}

If $n>k$, then $s\omega_k(P_n)=k$, hence the first inequality in Remark~\ref{rem:trivial} is sharp. If $n>k\ge 3$, then $s\omega_k(K_n-M)=n$, where $M$ is a matching of $K_n$. Hence the last two inequalities in Remark~\ref{rem:trivial} are also sharp. On the other hand, $\sgp_k(G)$ can be arbitrarily larger than $s\omega_k(G)$. For instance, if $r\ge 4$, then $s\omega_3(K_{1,r}) = 3$ and $\sgp_3(K_{1,r}) = r$.

The graphs attaining the equality in the rightmost inequality of Remark~\ref{rem:trivial} can be described as follows.

\begin{proposition}
\label{prop:charact-n-d-connected}
Let $G$ be a graph and let $k\in [2:n(G)-1]$. Then, $\sgp_k(G)=n(G)$ if and only if $G$ is $(n(G)-k+1)$-connected.
\end{proposition}

\proof
The statement $\sgp_k(G)=n(G)$ is equivalent to the fact that for every $k$-set $W$, a Steiner $W$-tree contains $k$ vertices. That every such subgraph $G[W]$ is connected is in turn equivalent to the fact that removing an arbitrary set of cardinality $n(G)-k$ does not disconnect the graph $G$, that is, $G$ is $(n(G)-k+1)$-connected.
\qed

Proposition~\ref{prop:charact-n-d-connected} in particular asserts that $\sgp_2(G)=n(G)$ if and only if $G$ is $(n(G)-1)$-connected, that is, if and only if $G$ is a complete graph. This fact was earlier observed in~\cite[Theorem 2.1]{ullas-2016}. By Proposition~\ref{prop:charact-n-d-connected}, this fact can be extended to the statement that the equality chain $\sgp_2(G)= \cdots =\sgp_{n(G)-1}(G)=n(G)$ holds if and only if $G$ is a complete graph. We next show that, a bit surprisingly, it is in general not true that a $k$-Steiner general position set is also a $k'$-Steiner general position set for some $k' > k$.

\begin{proposition}
\label{prop:k-1-and-not-k}
For every $k\ge 3$ there exist a graph $G^{(k)}$ and a set $A_k\subseteq V(G^{(k)})$ such that $A_k$ is a $(k-1)$-Steiner general position set and is not a $k$-Steiner general position set of $G^{(k)}$.
\end{proposition}

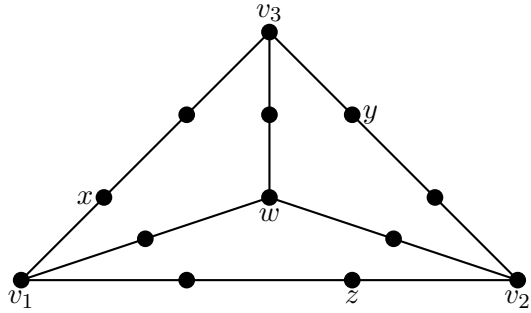
\begin{figure}[htb]
\begin{center}
\begin{tikzpicture}[scale=1.1,style=thick,x=1cm,y=1cm]
\def\vr{2.5pt} % \vr = vertex radius;

% define vertices
%%%%%
	%%% left
%%%%%
\path (0,0) coordinate (a);
\path (1,1) coordinate (x);
\path (2,2) coordinate (u1);
\path (3,3) coordinate (c);
\path (4,2) coordinate (y);
\path (5,1) coordinate (u2);
\path (6,0) coordinate (b);
\path (2,0) coordinate (u3);
\path (4,0) coordinate (z);

\path (3,1) coordinate (d);
\path (3,2) coordinate (u4);
\path (1.5,0.5) coordinate (u5);
\path (4.5,0.5) coordinate (u6);

%  edges
	\draw (a)--(c);
	\draw (a)--(d);
	\draw (c)--(d);
	\draw (b)--(d);
	\draw (a)--(b);
	\draw (b)--(c);

			\draw (a)[fill=black] circle (\vr);
			\draw (b)[fill=black] circle (\vr);
			\draw (c)[fill=black] circle (\vr);
			\draw (d)[fill=black] circle (\vr);
			\draw (x)[fill=black] circle (\vr);
			\draw (y)[fill=black] circle (\vr);
			\draw (z)[fill=black] circle (\vr);
			\draw (u1)[fill=black] circle (\vr);
			\draw (u2)[fill=black] circle (\vr);
			\draw (u3)[fill=black] circle (\vr);
			\draw (u4)[fill=black] circle (\vr);
			\draw (u5)[fill=black] circle (\vr);
			\draw (u6)[fill=black] circle (\vr);

\draw[anchor = north] (a) node {$v_1$};
\draw[anchor = north] (b) node {$v_2$};
\draw[anchor = south] (c) node {$v_3$};
\draw[anchor = north] (d) node {$w$};
\draw[anchor = east] (x) node {$x$};
\draw[anchor = west] (y) node {$y$};
\draw[anchor = north] (z) node {$z$};

\end{tikzpicture}
\end{center}
\caption{The graph  $G^{(3)}$. Also, the 2-sgp set $\{v_1,v_2,v_3,w\}$ is not 3-sgp set, and $\{w,x,y,z\}$ is a $k$-sgp set for $k\in \{2,3\}$.}\label{figure}
\end{figure}

\proof
Let $k\ge 3$ and construct $G^{(k)}$ as follows. First take the join of the cycle $C_k$ with consecutive vertices $v_1, v_2, \ldots, v_k$ and the one vertex graph $K_1$ with the vertex $w$. This creates a wheel $W_{k+1}$ with the center $w$, cf.\ Section~\ref{sec:joins}. Then subdivide each of the edges of the cycle $C_k$ by $k-1$ vertices and subdivide each of the spokes of the wheel with $k-2$ vertices. See Fig.~\ref{figure} where the graph $G^{(3)}$ is drawn.

We now claim that the set $A_k = \{w, v_1, v_2, \ldots, v_k\}$ is a $(k-1)$-Steiner general position set of $G^{(k)}$ and is not a $k$-Steiner general position set  of $G^{(k)}$.

Clearly, $|A_k| = k+1$. Let $B_k  = A_k \setminus \{w\}$, so that $|B_k| = k$. A smallest spanning tree that contains the vertices from $B_k$ and does not contain the vertex $w$ proceeds along the subdivided cycle $C_k$ and is of size $k(k-1)$. On the other hand, a spanning tree that contains the vertices from $B_k$ as well as the vertex $w$, contains all the subdivided spokes and is of size $k(k-1)$. Hence $d_{G^{(k)}}(B_k) = k(k-1)$ and the second described spanning tree implies that $A_k$ is not a $k$-Steiner general position set  of $G^{(k)}$.

Consider next $(k-1)$-subsets $B$ of $A_k$.  By symmetry, there are only two cases to consider. Suppose first that $w\in B$, so that $B$ contains $w$ and $k-2$ vertices of $C_k$. Then it is clear that the unique Steiner $B$-tree contains $w$ and the subdivided spokes between $w$ and the other vertices from $B$. (Its size is $(k-1)(k-2)$.) Suppose second that $B$ consists of $k-1$ vertices from $C_k$. Then a smallest spanning tree that contains vertices of $B$ and does not contain $w$ is of size $(k-2)k$. On the other hand, a spanning tree that contains $w$ and the subdivided spokes between $w$ and the vertices from $B$ is of size $(k-1)(k-1)$. Since $k(k-2) = k^2 - 2k < k^2 - 2k + 1 = (k-1)(k-1)$, we see that $d_{G^{(k)}}(B) = k(k-2)$ and conclude that $A_k$ is a $(k-1)$-Steiner general position set of $G^{(k)}$.
\qed

The proposition above asserts that there exist graphs containing $(k-1)$-sgp sets that are not $k$-Steiner general position sets for $k\ge 2$. However, there could yet exist other sets in such graphs that are $(k-1)$-sgp sets as well as $k$-Steiner general position sets, as for instance the set $\{w,x,y,z\}$ of Fig.~\ref{figure}, which is a $k$-sgp set for $k\in \{2,3\}$. In consequence, although there is no monotonicity with respect to inclusion for $k$-Steiner general position sets, it still could be a monotonicity relation with respect to the value of $\sgp_k(G)$ for every connected graph $G$. Proving or disproving a general monotonicity relation like this one seems to be a challenging problem.

The graph $G^{(3)}$ from Fig.~\ref{figure} was presented in~\cite{changat-2010} as a first in a family of graphs that has a 2-Steiner convex set that is not 3-Steiner convex. Later it was also mentioned in~\cite{ACKP}.

%%%%%%%%%%%%%%%%%%%%%%%%%%%%%%%%%%%%%%%%%%%%%%%%%%%%%
%%%%%%%%%%%%%%%%%%%%%%%%%%%%%%%%%%%%%%%%%%%%%%%%%%%%%
\section{Trees and cycles}
\label{sec:trees-cycles}
%%%%%%%%%%%%%%%%%%%%%%%%%%%%%%%%%%%%%%%%%%%%%%%%%%%%%
%%%%%%%%%%%%%%%%%%%%%%%%%%%%%%%%%%%%%%%%%%%%%%%%%%%%%

In this section we determine the $k$-Steiner general position number of trees and cycles.

\begin{theorem}
If $T$ is a tree with $n(T)\ge 3$ and $k\in [2:n(T)-1]$, then
$$\sgp_k(T)=\left\{\begin{array}{ll}
\ell(T); & k\leq \ell(T), \\[0.1cm]
k; & k>\ell(T).
\end{array}
\right.$$
\end{theorem}

\proof
Let $T$ be a tree of order at least $3$, and let $k\in [2:n(T)-1]$. Clearly, $L(T)$ is a $k$-Steiner general position set. From Remark~\ref{rem:trivial} we also know that $\sgp_k(T) \ge k$. Therefore,  $\sgp_k(T) \ge \max \{\ell(T), k \}$. To establish a corresponding lower bound, we distinguish two cases.

Assume first that $k\le \ell(T)$. Let us suppose that $\sgp_k(T) > \ell(T)$, and let $S$ be a $k$-sgp-set of $T$. Since $|S| > \ell(T)$, there exist three vertices $u,v,w\in S$ such that the shortest $u,v$-path in $T$ contains $w$. By taking a set $A\subseteq S$ of $k$ vertices, including $u,v$ and not including the vertex $w$, we obtain that the Steiner $A$-tree contains the vertex $w$, which is not possible. Hence, we conclude that $\sgp_k(T) \le \ell(T)$, and thus  $\sgp_k(T) = \ell(T)$, in the case when $k\le \ell(T)$.

Assume now that $k > \ell(T)$. By using a similar argument as in the previous paragraph, we find three vertices $u,v,w$ (and some other ones which could also exist) that allow to claim that there cannot be more than $k$ vertices in any $k$-sgp-set of $T$. For otherwise, we will obtain a Steiner tree with a not allowed vertex. Therefore, we deduce $\sgp_k(T)\le k$, which leads to the conclusion that if $k > \ell(T)$, then $\sgp_k(T) = k$.
\qed

\begin{theorem}
\label{thm:cycle}
If $n\ge 3$ and $k\in[2:n-1]$, then
$$\sgp_k(C_n)=\left\{\begin{array}{ll}
                       k; & k\in\left[\left\lfloor\frac{2n}{3}\right\rfloor :n-2\right], \\
                       k+1; & \mbox{otherwise}.
                     \end{array}
\right.$$
\end{theorem}

\proof
Let $C_n=v_0\cdots v_{n-1}v_0$. In the rest of the proof all the operations with the indices of the vertices of $C_n$ are done modulo $n$.

Consider a set $A\subseteq V(C_n)$ with $|A| = k+2$. Let $a,b,c$ be three consecutive vertices of $A$, where consecutive refers to their order on $C_n$. Then $B=A-\{a,c\}$ is a $k$-set. Since every Steiner $B$-tree contains at least one of the vertices $a$ or $c$, we infer that $A$ is not a $k$-Steiner general position set. By Remark~\ref{rem:trivial} we have $\sgp_k(C_n)\ge k$, hence $\sgp_k(C_n) \in \{k, k+1\}$ follows.

As $C_n$ is $2$-connected, Proposition \ref{prop:charact-n-d-connected} implies that $\sgp_{n-1}(C_n)=n$. In the rest we may thus assume that $k\le n-2$.

We claim that if $S$ is a $k$-Steiner general position set of $C_n$ of cardinality $k+1$, then $S$ contains no three consecutive vertices of $C_n$. Suppose on the contrary that, without loss of generality, $S$ contains $v_0$, $v_1$, and $v_2$. Then for the set $S'=S\setminus\{v_1\}$, there is at least one Steiner $S'$-tree which contains the vertex $v_1$, a contradiction proving the claim. From it and by the pigeonhole principle, if $\sgp_k(C_n)=k+1$, then $k+1\le \left\lfloor\frac{2n}{3}\right\rfloor$. We have thus proved that $\sgp_k(C_n)=k$ for every $k\in \left[\left\lfloor\frac{2n}{3}\right\rfloor :n-2\right]$.

In the rest of the proof let $k\in [2:\left\lfloor\frac{2n}{3}\right\rfloor-1]$. We want to show that in these cases $\sgp_k(C_n) = k+1$. Since we already know that $\sgp_k(C_n) \le k+1$ it remains to construct a $k$-Steiner general position set of cardinality $k+1$. For this sake write $n$ as $n=q(k+1)+r$, where $r< k+1$, and consider the following two cases.

Assume first that $r<q$. Let
$$A=\{v_0,v_{q},v_{2q},\dots,v_{kq}\}$$
and note that $|A| = k+1$. If $v_{i_1},v_{i_2},v_{i_3}$ are any three consecutive vertices of $A$, then we note that $d_{C_n}(v_{i_1},v_{i_3})\in \{2q, 2q+r\}$, which is strictly larger than $d_{C_n}(v_{j},v_{l})\in \{q, q+r\}$, where $v_j$ and $v_l$ are arbitrary consecutive vertices of $A$. In consequence, we deduce that for any vertex $v_j\in A$, the Steiner $(A\setminus\{v_j\})$-tree does not contain the vertex $v_j$. Thus, $A$ is a $k$-Steiner general position set of $C_n$ and so, $\sgp_k(C_n)\ge k+1$ as required.

Assume now that $r\ge q$. The idea is then to make a partition of $V(C_n)$ into $k+2$ sets of consecutive vertices starting with $v_0$: $r-q+1$ sets of cardinality $q+1$; $k+1-r+q-1$ sets of cardinality $q$; and one set of cardinality $q-1$. Take first  vertex from each of the first $k+1$ sets (and no vertex from the set of cardinality $q-1$). The set obtained is:
$$A'=\{v_0,v_{q+1},\dots , v_{(r-q+1)(q+1)}, v_{(r-q+1)(q+1)+q}, \dots, v_{(r-q+1)(q+1)+(k+1-r+q-2)q}\}\,.$$
Notice that the last element of $A'$ is $v_{(r-q+1)(q+1)+(k+1-r+q-2)q}=v_{n-2q+1}$. If $v_{i_1}$, $v_{i_2}$, and $v_{i_3}$ are arbitrary three consecutive vertices of $A'$, then $d_{C_n}(v_{i_1},v_{i_3})\in \{2q, 2q+1, 2q+2, 3q\}$. Moreover, if $v_j$ and $v_l$ are two consecutive vertices of $A'$, then $d_{C_n}(v_{j},v_{l})\in \{q, q+1, 2q-1\}$. Thus, $d_{C_n}(v_{i_1},v_{i_3})$ is always strictly larger than $d_{C_n}(v_{j},v_{l})$. As a consequence, we again deduce that for any vertex $v_j\in A$, the Steiner $(A\setminus\{v_j\})$-tree does not contain the vertex $v_j$. Theorefore, $A$ is a $k$-Steiner general position. Hence, also in this second case we have $\sgp_k(C_n)\ge k+1$ and we are done.
\qed

%%%%%%%%%%%%%%%%%%%%%%%%%%%%%%%%%%%%%%%%%%%%%%%%%%%%%
%%%%%%%%%%%%%%%%%%%%%%%%%%%%%%%%%%%%%%%%%%%%%%%%%%%%%
\section{Joins of graphs}
\label{sec:joins}
%%%%%%%%%%%%%%%%%%%%%%%%%%%%%%%%%%%%%%%%%%%%%%%%%%%%%
%%%%%%%%%%%%%%%%%%%%%%%%%%%%%%%%%%%%%%%%%%%%%%%%%%%%%

In this section we give a formula for the $k$-Steiner general position number of graph joins. For the Steiner diameter of joins see~\cite{wang-2019}, and for the Steiner diameter of several other graph operations~\cite{mao-2017, wang-2019}.

The join  $G\vee H$ of disjoint graphs $G$ and $H$ is the graph with vertex set $V(G\vee H)=V(G)\cup V(H)$ and edge set
$$E(G\vee H)=E(G)\cup E(H)\cup \{gh:g\in V(G),h\in V(H)\}.$$
Known families of graphs that can be presented as the join of two graphs include complete graphs $K_n=K_p\vee K_{n-p}$, complete bipartite graphs $K_{s,t}=\overline{K}_s\vee \overline{K}_t$, wheel graphs $W_r=K_1\vee C_{r-1}$, $r\geq 4$, and fan graphs $F_n=K_1\vee P_{n-1}$, $n\geq 2$.

It is well-known that $\omega(G\vee H) =\omega(G) + \omega(H)$ and also not difficult to extend this result to $s\omega_k(G\vee H)=s\omega_k(G)+s\omega_k(H)$ for connected graphs $G$ and $H$. In the general case the corresponding result still holds as shown next.

\begin{lemma}
\label{lem:join}
If $G$ and $H$ are graphs and $k\in [2:n(G)+n(H)]$, then
$$s\omega_k(G\vee H)=s\omega_k(G)+s\omega_k(H).$$
\end{lemma}

\proof
If $A\subseteq V(G)$ and $B\subseteq V(H)$ are $k$-Steiner cliques of $G$ and $H$, respectively, then $A\cup B$ is a $k$-Steiner clique of $G\vee H$ even if $G$ or $H$ are disconnected and every component has less than $k$ vertices. In this latter case, one selects any arbitrary $\min\{n(G),k-1\}$ vertices in $G$, which together with a $k$-Steiner clique in $H$ form a $k$-Steiner clique of $G\vee H$. A similar argument can be used symmetrically for $H$. Thus, $s\omega_k(G\vee H)\ge  s\omega_k(G)+s\omega_k(H)$.

Suppose that $s\omega_k(G\vee H)>s\omega_k(G)+s\omega_k(H)$ and let $S$ be $k$-Steiner clique of $G\vee H$. Hence, it must happen that $|S\cap V(G)|>s\omega_k(G)$ or $|S\cap V(H)|>s\omega_k(H)$. Without loss of generality, we may assume that $|S\cap V(G)|>s\omega_k(G)$. Then $S\cap V(G)$ is not a $k$-Steiner clique of $G$, and there exists a $k$-subset $B$ of $S\cap V(G)$, where $G[B]$ is not connected. This also means $(G\vee H)[B]$ is not connected, and $S$ is not a $k$-Steiner clique of $G\vee H$, a contradiction. Hence, $s\omega_k(G\vee H)\leq s\omega_k(G)+s\omega_k(H)$ and the equality follows.
\qed

Let $G$ be a graph and let $A\subseteq V(G)$ be a set of vertices of cardinality at least $k$. Then we say that $A$ is a $k$-\textit{Steiner join-critical set} of $G$ if for each $k$-subset $B$ of $A$ we have $d_{G[A]}(B)\neq k$. That is, $A$ is a $k$-Steiner join-critical set if there exists no $k$-set $B\subseteq A$ such that a Steiner $B$-tree in $G[A]$ contains $k+1$ vertices. By $\sjc_k(G)$ we denote the cardinality of a largest $k$-Steiner join-critical set. For a given set $D\subseteq V(G)$, if every connected component of $G[D]$ is of order at most $k$, then $D$ is $k$-Steiner join-critical. For the particular case in which $D=V(G)$, if every connected component of $G$ has order at most $k$, then $\sjc_k(G) = n(G)$. Note also that, by definition, if $k \ge n(G)$, then $\sjc_k(G) = n(G)$.

It seems that determining $\sjc_k(G\vee H)$ is a hard problem, but one can express the exact value for $\sgp_k(G\vee H)$ in terms of $\sjc_k(G)$, $\sjc_k(H)$, and $s\omega_k(G\vee H)$. In addition, by Lemma~\ref{lem:join}, the latter invariant can also be expressed by related invariants of $G$ and $H$.

\begin{theorem}\label{join}
If $G$ and $H$ are graphs and $k\in[2:n(G\vee H)-1]$, then
$$\sgp_k(G\vee H)=\max\{s\omega_k(G\vee H),\sjc_k(G), \sjc_k(H)\}.$$
\end{theorem}

\proof
Let $G$ and $H$ be graphs and let $M=\max\{s\omega_k(G\vee H),\sjc_k(G), \sjc_k(H)\}$. If $M=s\omega_k(G\vee H)$, then $\sgp_k(G\vee H)\geq M$ by Remark \ref{rem:trivial}. Suppose next that $M=\sjc_k(G)$. Let $A\subseteq V(G)$ be a $k$-Steiner join-critical set of $G$ of cardinality $\sjc_k(G)$. We wish to show that $A$ is a $k$-Steiner general position set of $G\vee H$. Let $B$ be any $k$-subset of $A$. If $G[B]$ is connected, then we are done. Hence assume that $G[B]$ is not connected. Let $h$ be an arbitrary vertex of $H$. Then the set $B_h=B\cup\{h\}$ induces a connected subgraph of $G\vee H$  which means that $d_{G\vee H}(B)\leq k$. On the other hand, there does not exists a Steiner $B$-tree of order $k+1$ that contains only vertices of $A$, because $A$ is $k$-Steiner join-critical in $G$. So, every Steiner $B$-tree in $G\vee H$ does not contain any additional vertex from $A$, and consequently $A$ is a $k$-Steiner general position set for $G\vee H$, meaning that $\sgp_k(G\vee H)\geq M$ also holds in this case. By the symmetry of $G$ and $H$ in $G\vee H$, we also get that $\sgp_k(G\vee H)\geq M$ when $M=\sjc_k(H)$.

Suppose now that there exist a $k$-Steiner general position set $A$ of $G\vee H$ of cardinality greater than $M$. Assume first that $A_H=A\cap V(H)\neq \emptyset$ and $A_G=A\cap V(G)\neq \emptyset$. Since $|A|>M\geq s\omega_k(G\vee H)$, Lemma \ref{lem:join} implies that at least one of $A_G$ and $A_H$ contains more vertices than $s\omega_k(G)$ and $s\omega_k(H)$, respectively. Assume without loss of generality that $|A_G|>s\omega_k(G)$. Recall that for any graph $G$ either
\begin{itemize}
  \item $s\omega_k(G)=\min\{n(G), k-1\}$ (which means that every connected component of $G$ has cardinality smaller than $k$), or
  \item $s\omega_k(G)\geq k$.
\end{itemize}
Let $M_1=\min\{n(G), k-1\}$. If $M_1=k-1$, then in both situations above, $|A_G|>M_1$ implies that $|A_G|\geq k$. This means there exists a $k$-subset $B$ of $A_G$, such that $G[B]$ is not a connected subgraph and $d_{G\vee H}(B)\geq k$. The set $B_h=B\cup\{h\}$, with $h\in A_H$, induces a connected subgraph of $G\vee H$ and $d_{G\vee H}(B_h)=k$. So, $B_h$ induces a Steiner $B$-tree that contains a vertex from $A\setminus B$, a contradiction with $B$ being a $k$-Steiner general position set of $G\vee H$. So, let now $M_1 = n(G) < k$. In this case we have $|A_G|> M_1 = n(G)$, a contradiction again.

It remains to consider the case when $A_G=\emptyset$ or $A_H=\emptyset$. We may, without loss of generality, assume that $A_H=\emptyset$. Because $|A|>M\geq \sjc_k(G)$, $A$ is not $k$-Steiner join-critical and there exists a $k$-subset $B$ of $A$, for which $d_{G[A]}(B)=k$. Clearly, $G[B]$ is not connected and $d_{G\vee H}(B)=k$, a contradiction with $A$ being a $k$-Steiner general position set. Hence $\sgp_k(G\vee H)\leq M$ and the equality follows.
\qed

To give some applications of Theorem~\ref{join}, we first determine exact results for $\sjc_k(P_n)$ and $\sjc_k(C_n)$.

\begin{proposition}\label{pathcycle}
Let $n\ge 3$. If $k\in[2:n-1]$, then  $\sjc_k(P_n)=n-\left\lfloor \frac{n}{k+1}\right\rfloor$. If $\ell\in[2:n-2]$, then  $\sjc_{\ell}(C_n)=n-1-\left\lfloor \frac{n-1}{\ell+1}\right\rfloor$. Moreover, $\sjc_{n-1}(C_n)=n$.
\end{proposition}

\proof
Divide the vertex set of $P_n$ into $\left\lfloor \frac{n}{k+1}\right\rfloor$ sets of $k+1$ consecutive vertices and a remainder set with at most $k$ vertices. Let $Q$ be the set consisting of the last vertex from each of the $\left\lfloor \frac{n}{k+1}\right\rfloor$ sets of $k+1$ consecutive vertices. Every connected component of the subgraph induced by $B=V(P_n)-Q$ has at most $k$ vertices, therefore $B$ is a $k$-Steiner join-critical set of cardinality $n-\left\lfloor \frac{n}{k+1}\right\rfloor$. If $\sjc_k(P_n)>n-\left\lfloor \frac{n}{k+1}\right\rfloor$ would hold, then every $k$-Steiner join-critical set $B$ of cardinality $\sjc_k(P_n)$ would contain a connected component $C$ with at least $k+1$ vertices. If $x$ is a middle vertex of such a component $C$, then $d_{P_n[B]}(C-\{x\})=k$, a contradiction. So, the equality holds for  paths.

For cycles we can use the same steps for $\ell\in[2:n-2]$, only that we need to put also the last vertex into $Q$, and then the result follows. If $\ell=n-1$, then any $\ell$ vertices of $C_n$ induce a connected graph and we are done.
\qed

We can now apply Theorem~\ref{join} to specific families of graphs as follows.

\begin{corollary}\label{exact} The following assertions hold for positive integers $k,n,r,s$.
\begin{itemize}
\item[(i)] If $k\in[2:r+s-1]$, then $\sgp_k(K_{r+s})=\sgp_k(K_r\vee K_s)=r+s$.
\item[(ii)] If $n\geq 6$ and $k\in[2:n-1]$, then $\sgp_k(W_n)=\sgp_k(K_1\vee C_{n-1})=\max\{k+1,n-2-\left\lfloor \frac{n-2}{k+1}\right\rfloor$\}.
\item[(iii)] If $n\geq 4$ and $k\in[2:n-1]$, then $\sgp_k(F_n)=\sgp_k(K_1\vee P_{n-1})=\max\{k+1,n-1-\left\lfloor \frac{n-1}{k+1}\right\rfloor$\}.
\item[(iv)] If $r\leq s$ and $k\in[2:r+s-1]$, then
 $$\sgp_k(K_{r,s})=\sgp_k(\overline{K}_r\vee \overline{K}_s)=\left\{\begin{array}{ll}
                                                                      \max\{s,\min\{k-1,r\}+k-1\}; & k\leq s, \\[0.1cm]
                                                                      r+s; & k>s.
                                                                    \end{array}
\right.$$
\item[(v)]  If $k\in[2:r+s-1]$, then
$$\sgp_k(K_r\vee \overline{K}_s)=\left\{\begin{array}{ll}
                                                                      r+\min\{s,k-1\}; & k>\min\{r,s\}, \\[0.1cm]
			\max\{r+k-1,s\}; & k\leq\min\{r,s\}.
                                                                    \end{array}
\right.$$
\end{itemize}
\end{corollary}

\proof
The results are obtained by straightforward applications  of Lemma~\ref{lem:join} and Theorem~\ref{join}. For items $(ii)$ and $(iii)$, Proposition~\ref{pathcycle} is also needed. For some cases, Proposition \ref{prop:charact-n-d-connected} can also be applied.
\qed

%%%%%%%%%%%%%%%%%%%%%%%%%%%%%%%%%%%%%%%%%%%%%%%%%%%%%
%%%%%%%%%%%%%%%%%%%%%%%%%%%%%%%%%%%%%%%%%%%%%%%%%%%%%
\section{Lexicographic products}
\label{sec:lexico}
%%%%%%%%%%%%%%%%%%%%%%%%%%%%%%%%%%%%%%%%%%%%%%%%%%%%%
%%%%%%%%%%%%%%%%%%%%%%%%%%%%%%%%%%%%%%%%%%%%%%%%%%%%%

Let $G$ and $H$ be two graphs. The lexicographic product $G\circ H$ is a graph with $V(G\circ H)=V(G)\times V(H)$. Two vertices $(g,h)$ and $(g',h')$ are adjacent if $gg'\in E(G)$ or ($g=g'$ and $hh'\in E(H$)). The lexicographic product is a kind of generalization of join because $K_2\circ G\cong G\vee G$ for any graph $G$. The map $p_G:(g,h)\mapsto g$ is the projection of $V(G\circ H)$ to $V(G)$. The set $G^h=\{(g,h):g\in V(G)\}$ is called the $G$-\emph{layer} (through $h$). Similarly, $^gH=\{(g,h):h\in V(H\}$ is called the $H$-\emph{layer} (through $g$). The subgraphs of $G\circ H$ induced by $G^h$ and by $^gH$ are clearly isomorphic to $G$ and $H$, respectively.

Steiner trees are behaving relatively nice with respect to the first factor $G$ of lexicographic product. More accurate, the following lemma was proved in \cite[Lemma 3.1]{ACKP}.

\begin{lemma}
\label{simple}Let $g_{1},\ldots ,g_{k}$ be different vertices of a
connected graph $G$. Then for any (not necessarily different) vertices $%
h_{1},,\ldots ,h_{k}$ of a graph $H$, a Steiner tree of $%
g_{1},\ldots ,g_{k}$ (in $G$) and a Steiner tree of $%
(g_{1},h_{1}),\ldots ,(g_{k},h_{k})$ (in $G\circ H$) have the
same size.
\end{lemma}

This forms a basis for any set $B=\{(g_{1},h_{1}),\ldots ,(g_{k},h_{k})\}$ of those vertices that project to at least two different vertices of $G$. Namely, we have the same size as the Steiner tree of $p_G(B)=\{g_1,\ldots, g_k\}$ plus $m_i-1$, for every $i\in [k]$, where $m_i$ represents the number of vertices from $B$ that project to $g_i$. So, the only case that is not connected with the Steiner tree in $G$ occurs when all the vertices of $B$ project to the same vertex $g_1$.

Let $k$ and $\ell $ be two positive integers with $k\leq\ell$. A set $A\subseteq V(G)$ is a $[k:\ell]$-\emph{Steiner general position set} of a graph $G$, or a $[k:\ell]$-sgp set for short, if it is a $j$-Steiner general position set for every $j\in[k:\ell]$. The cardinality of a largest  $[k:\ell]$-sgp set for $G$ is represented as $\sgp_{[k:\ell]}(G)$. The family of all $[k:\ell]$-sgp sets of a graph $G$ is denoted by ${\cal G}_{k,\ell}$. For every $k$-general Steiner position set $S$, we partition $S$ into two sets ${\cal I}_S$ and ${\cal J}_S$, where ${\cal I}_S$ contains all isolated vertices in the subgraph of $G[S]$ and ${\cal J}_S=S\setminus {\cal I}_S$. Every set of vertices of cardinality at most $k$ is a $[k:\ell]$-sgp set. On the other hand, as can be seen from the proof of Theorem~\ref{thm:cycle}, for a cycle $C_n$, any set with $k+2$ vertices is not a $k$-sgp set. So, any $[k:\ell]$-sgp set of $C_n$ contains at most $k+1$ vertices. Fig.~\ref{figure} shows a graph where the set $S_1=\{v_1,v_2,v_3,w\}$ is a 2-sgp set, but not a 3-sgp set. However, the set $S_2=\{x,y,z,w\}$ is a $[2:13]$-sgp set. Clearly, ${\cal I}_{S_1}=S_1$ and ${\cal I}_{S_2}=S_2$, while ${\cal J}_{S_1}=\emptyset={\cal J}_{S_2}$.

\begin{theorem}\label{lex}
Let $G$ and $H$ be nontrivial graphs, let $G$ be connected, and let $k\in[2:n(G)\cdot n(H)-1]$.
If $j=\left\lceil \frac{k}{n(H)}\right\rceil$ and $\ell=\min\{k,n(G)\}$, then
$$\sgp_k(G\circ H)\geq\left\{\begin{array}{ll}
                                                \max_{S\in {\cal G}_{2,k}}\{|{\cal I}_S|\sjc_k(H)+|{\cal J}_S|s\omega_k(H)\}; &  k\leq n(H), \\[0.1cm]
                                                 \sgp_{[j:\ell]}(G)n(H); & n(H)<k<n(G)\cdot n(H).
                                                                    \end{array}
\right.$$
Moreover, the equality holds if $k>(n(G)-1)n(H)$.
\end{theorem}

\proof
Let first $k\leq n(H)$ and let $M=\max_{S\in {\cal G}_{2,k}}\{|{\cal I}_S|\sjc_k(H)+|{\cal J}_S|s\omega_k(H)\}$. Fix an arbitrary $S\in {\cal G}_{2,k}$ together with ${\cal I}_S$ and ${\cal J}_S$. In addition, let $D_1$ be a $\sjc_k(H)$-set and let $D_2$ be an $s\omega_k(H)$-set. We will show that $A=({\cal I}_S\times D_1)\cup ({\cal J}_S\times D_2)$ is a $k$-sgp set of $G\circ H$. Let $B$ be any $k$-subset of $A$. Suppose first that all vertices of $B$ belong to one layer $^gH$. Clearly, any Steiner $B$-tree contains either $k$ or $k+1$ vertices. Moreover, we have $d_H(B)=k-1$ when $p_H(B)$ induces a connected subgraph. If $g\in {\cal J}_S$, then $B$ induces a connected subgraph since $D_2$ is an $s\omega_k(G)$-set and we are done. If $g\in {\cal I}_S$, then again we are done when $d_H(p_H(B))=k-1$, since $B$ induces a connected graph. Otherwise, let $d_H(p_H(B))\geq k+1$ (recall $D_1$ is a $\sjc_k(H)$-set), and the additional vertex of a Steiner $B$-tree, say $(x,y)$, belongs either to $^gH$ or to $^{g'}H$ for some neighbor $g'$ of $g$. If $x=g$, then $y\notin D_1$ by the definition of a $\sjc_k(H)$-set and $(x,y)\notin A$. If $x=g'$, then $(x,y)\notin A$ because $g\in {\cal I}_S$. Thus, the only vertices from $A$ that are included at some Steiner $B$-tree are only those ones from $B$.

Suppose now that $|p_G(B)|=j>1$, that is, $j\in [2:k]$. As $S\in {\cal G}_{2,k}$, any Steiner $p_G(B)$-tree in $G$ contains no other vertices from $S$ than those ones from $p_G(B)$. By a consequence of Lemma \ref{simple}, we get that every Steiner $B$-tree in $G\circ H$ contains no other vertices from $A$ than those ones from $B$. Hence, $A$ is a $k$-sgp set and we have $\sgp_k(G\circ H)\geq M$.

%To show the equality let $A$ be a $\sgp_k(G\circ H)$-set and let $A_G=p_G(A)$. We denote by ${\cal I}_{A_G}$ all isolated vertices from $G[A_G]$ and let ${\cal J}_{A_G}=A_G\setminus{\cal I}_{A_G}$. Suppose first that there exists $g\in {\cal I}_{A_G}$ such that $|^gH\cap A|>\sjc_k(H)\geq k$. By the definition of $\sjc_k(H)$ there exists $k$-subset $B$ of $^gH\cap A$, such that a Steiner $B$-tree contains a vertex from $(^gH\cap A)\setminus B$, a contradiction with $A$ being a $\sgp_k(G\circ H)$-set. So, $|^gH\cap A|\leq\sjc_k(H)$ for every $g\in {\cal I}_{A_G}$.

%Assume now that there exists $g\in {\cal J}_{A_G}$ such that $|^gH\cap A|>s\omega_k(H)\geq k-1$. By the definition of $s\omega_k(H)$, there exists a $k$-subset $B$ of $^gH\cap A$, such that $p_H(B)$ is not connected in $H$. Therefore, $B\cup\{(g',h)\}$ forms a Steiner $B$-tree where $g'\in {\cal J}_{A_G}$ is a neighbor of $g$, and $(g',h)\in A$, the same contradiction again. Hence, for every $\sgp_k(G\circ H)$-set $A$, we have $|A|\leq|{\cal I}_{A_G}|\sjc_k(H)+|{\cal J}_{A_G}|s\omega_k(H)$. Moreover, if $|A|<|{\cal I}_{A_G}|\sjc_k(H)+|{\cal J}_{A_G}|s\omega_k(H)$, then we have a contradiction since every

%We still need to show that $A_G\in {\cal G}_{2,k}$. If not, then there exists $j\in [2:k]$ and a $j$-subset $B_G=\{g_1,\cdots,g_j\}$ of $A_G$, such that a Steiner tree for $B_G$ contains some at least one vertex $g_0$ from $A_G\setminus B_G$. Hence there exists SIMMILAR PROBLEM AS TWO PARAGRAPHS BELOW. SEE COMMENT THERE

Let now $n(H)<k<n(G)\cdot n(H)$, $j=\left\lceil \frac{k}{n(H)}\right\rceil$, and $\ell=\min\{k,n(G)\}$. We will show that $\sgp_k(G\circ H)\geq \sgp_{[j:\ell]}(G)n(H)$ by proving that $A=A_G\times V(H)$ is a $k$-Steiner general position set of $G\circ H$, where $A_G$ is a $\sgp_{[j:\ell]}(G)$-set. Let $B$ be any $k$-subset of $A$ and let $B_G=p_G(B)$. Clearly, $j\leq|B_G|\leq \ell$ and every Steiner $B_G$-tree does not contain any additional vertex from $A_G$ because $A_G$ is a $\sgp_{[j:\ell]}(G)$-set. But then by Lemma \ref{simple} and the comment after it, also the Steiner $B$-tree does not contain any additional vertex from $A$. Thus, $A$ is a $k$-Steiner general position set of $G\circ H$ and the first inequality follows.

If $k>(n(G)-1)n(H)$, then $p_G(B) = V(G)$ for any $k$-set $B$. Since $G$ is connected, $(G\circ H)[B]$ is also connected and every $B$-Steiner tree contains only vertices from $B$. Therefore, the equality holds.
\qed

%Suppose now that $\sgp_k(G\circ H)> \sgp_{[j:\ell]}(G)n(H)$. Let $A$ be a $\sgp_k(G\circ H)$-set and let $A_G=p_G(A)$. By our assumption $|A_G|>\sgp_{[j:\ell]}(G)$, and there exists $t\in[j:\ell]$ such that $A_G$ is not a $t$-general Steiner position set. Hence, there exists a $t$-subset $B_G$ of $A_G$ such that some Steiner $B_G$-tree contains a vertex $g$ from $A_G\setminus B_G$. We denote by $(g,h)$ a vertex from $A\cap ^gH$. If there exists at least $k$ vertices in $A\cap (B_G\times V(H))$, then we can choose a $k$-subset $B$ with $p_G(B)=B_G$. By Lemma \ref{simple}, and the comment after it, there exists a Steiner $B$-tree that contains $(g,h)\in A\setminus B$, a contradiction with $A$ being a $\sgp_k(G\circ H)$-set. Therefore $A\cap (B_G\times V(H))$ contains less than $k$ vertices. PROBLEMS!!! I HAVE NO IDEA WHY THIS LEADS TO A CONTRADICTION. IT SEEMS THAT IT COULD BE POSSIBLE TO HAVE SMALL AMOUNT OF VERTICES IN EVERY H LAYER WHEN THE ORDER OF H IS SMALL AND THAT THIS BRINGS MORE THAN IF YOU TAKE WHOLE H LAYERS. HOWEVER I BELIEVE THAT EQUALITY HOLDS, BUT DO NOT KNOW HOW TO PROVE IT.

If we set $k=2$, then we can show the equality in Theorem \ref{lex}. For this, notice that in any 2-Steiner join-critical set $A$, two nonadjacent vertices cannot have any of their common neighbors also in $A$. Hence, $G[A]$ is a disjoint union of complete graphs. Moreover, $S\in{\cal G}_{2,2}$ simply means that $S$ is a general position set of $G$. The study of the general position number of the lexicographic product of graphs was initiated in \cite{klavzar-2019}, but just finding a connection between such parameter and other related structure. We next give a formula for the general position number of this product.

\begin{theorem}\label{lexgp}
If $G$ and $H$ are nontrivial graphs where $G$ is connected, then
$$\gp(G\circ H)=\max_{S\in {\cal G}_{2,2}}\{|{\cal I}_S|\sjc_2(H)+|{\cal J}_S|\omega(H)\}.$$
\end{theorem}

\proof
By Theorem \ref{lex} we have  $\gp(G\circ H)\geq \max_{S\in {\cal G}_{2,2}}\{|{\cal I}_S|\sjc_2(H)+|{\cal J}_S|\omega(H)\}$.

To show the equality, let $A$ be a $\gp(G\circ H)$-set and let $A_G=p_G(A)$. We denote by ${\cal I}_{A_G}$ the set of isolated vertices from $G[A_G]$, and let ${\cal J}_{A_G}=A_G\setminus{\cal I}_{A_G}$. Suppose first that there exists $g\in {\cal I}_{A_G}$ such that $|^gH\cap A|>\sjc_2(H)\geq 2$. By the definition of $\sjc_2(H)$, there exists a $2$-subset $B$ of $^gH\cap A$, such that any Steiner $B$-tree contains a vertex from $(^gH\cap A)-B$, a contradiction with $A$ being a $\gp(G\circ H)$-set. Thus, $|^gH\cap A|\leq\sjc_2(H)$ for every $g\in {\cal I}_{A_G}$.

Assume now that there exists $g\in {\cal J}_{A_G}$ such that $|^gH\cap A|>\omega(H)\geq 1$. By the definition of $\omega(H)$, there exists a $2$-subset $B$ of $^gH\cap A$, such that $p_H(B)$ is not connected in $H$. Therefore, $B\cup\{(g',h)\}$ forms a Steiner $B$-tree, where $g'$ is a neighbor of $g$ in ${\cal J}_{A_G}$ and $(g',h)\in A$, the same contradiction again. Hence, for every $\gp(G\circ H)$-set $A$ we have $|A|\leq|{\cal I}_{A_G}|\sjc_2(H)+|{\cal J}_{A_G}|\omega(H)$.

We still need to show that $A_G\in {\cal G}_{2,2}$. If not, then there exists a $2$-subset $B_G=\{g,g'\}$ of $A_G$ such that a $g,g'$-geodesic contains another vertex $g_0$ from $A_G$. Let $h,h',h_0\in V(H)$ be such that $(g,h),(g',h'),(g_0,h_0)\in A$. Clearly $(g_0,h_0)$ belongs to a $(g,h),(g',h')$-geodesic in $G\circ H$, a contradiction with $A$ being a $\gp(G\circ H)$-set.
\qed

Notice that the first two paragraphs of the proof above also suit for every $2< k\leq n(H)$. Unfortunately, there are several problems with the last paragraph while trying to get the equality in Theorem \ref{lex}.

We next illustrate that we need sets from ${\cal G}_{2,k}$ when $2< k\leq n(H)$. For this, notice that, by Theorem~\ref{thm:cycle}, the cardinality of a set $S$ such that $S\in {\cal C_n}_{2,k}$ is at most three for $n\geq 5$. By Theorem~\ref{lex}, we have
$$\sgp_k(C_n\circ K_\ell)\geq 3\ell.$$
On the other hand, if a general position set $S$ of $C_n\circ K_\ell$ projects to more than three vertices of $G$, then let $g_1,\dots,g_t$ be consecutive vertices of $p_G(S)$ on $C_n$. By $S_{i,j}$ we denote the subset of $S$ that projects to $P_G(S)-\{g_i,g_j\}$ where $|i-j|>1$ and $\{i,j\}\neq\{1,t\}$. If $|S_{i,j}|\geq k$, then there exists a $k$-subset $B$ of $S_{i,j}$ that projects to some vertices $g_p,g_r$ where $i<p<j$ and ($r<i$ or $r>j$). Clearly, there exists a $B$-Steiner tree that contains a vertex from $S\cap ^{g_i}H$ or from $S\cap ^{g_j}H$, a contradiction. So, $|S_{i,j}|<k$ and consequently $|S|<2k-1$. Hence $S$ is not a $\sgp_k(C_n\circ K_\ell)$-set. With this we have also shown the following.

\begin{corollary}\label{lexCK}
If $n\geq 5$ and $\ell\geq 2$, then
$$\sgp_k(C_n\circ K_\ell)=3\ell.$$
\end{corollary}

%%%%%%%%%%%%%%%%%%%%%%%%%%%%%%%%%%%%%%%%%%%%%%%%%%%%%
%%%%%%%%%%%%%%%%%%%%%%%%%%%%%%%%%%%%%%%%%%%%%%%%%%%%%
\section{Split graphs and infinite grids}
\label{sec:grids}
%%%%%%%%%%%%%%%%%%%%%%%%%%%%%%%%%%%%%%%%%%%%%%%%%%%%%
%%%%%%%%%%%%%%%%%%%%%%%%%%%%%%%%%%%%%%%%%%%%%%%%%%%%%

A connected graph $G$ is a \emph{split graph} if we can partition $V(G)$ in two sets $Q$ and $I$, such that $G[Q]$ is a clique $K_r$ and $G[I]$ is a graph without edges $\overline{K}_s$. By $G_{r,s}$ we denote a split graph with clique $K_r$ and independent set $\overline{K}_s$. Among joins, $K_r\vee \overline{K}_s$ is a split graph (see Corollary \ref{exact} $(v)$ for its $k$-Steiner general position number). We order the vertices of $I=\{v_1,\ldots,v_s\}$ by their degree, that is $\delta(v_1)\geq\cdots\geq\delta(v_s)$. Let $i(G)$ be defined as follows. If $\delta(v_1)\le r-k+1$, then $i(G)=0$. Otherwise, $i(G)$ denotes the largest integer $i\in[s]$ such that $\delta(v_i)>r-k+i$.

A \emph{universal vertex} of a graph $G$ is a vertex of degree $n(G)-1$. The set of all universal vertices in a graph $G_{r,s}$ is denoted by $U$. We also define
$$u_k(G_{r,s})=\left\{\begin{array}{lll}
                                                |U|; &  k>s, \\[0.1cm]
                                                 0;  &  k\leq s.
                                                                    \end{array}
\right.$$

\begin{theorem}\label{split}
If $G_{r,s}$ is a split graph, then
$$\sgp_k(G_{r,s})\geq\max\{r+i(G_{r,s}),s+u_k(G_{r,s}),k\}.$$
\end{theorem}

\proof
Let $G_{r,s}$ be a split graph with clique $Q$ on $r$ vertices and independent set $I$ on $s$ vertices and let $M=\max\{r+i(G_{r,s}),s+u_k(G_{r,s}),k\}$. If $M=k$, then the inequality holds by Remark \ref{rem:trivial} and we are done. 

Suppose next that $M=s+u_k(G_{r,s})$. This means that $s+u_k(G_{r,s})>k$, for otherwise we are in the case $M=k$ above. Add to the set $A$ all vertices from $I$ and all universal vertices when $k>s$. Select any $k$-subset $B$ of $A$. If $u_k(G_{r,s})=0$, then $A\subset I$ and every Steiner $B$-tree contains only vertices from $Q$, beside those already in $B$. Moreover, if there is a universal vertex in $A$, then there is at least one universal vertex also in $B$ because $k>s$. In this case, $B$ induces a connected subgraph and any Steiner $B$-tree does not contain any vertex from $A\setminus B$. Thus, $A$ is a $k$-Steiner general position set for $G_{r,s}$ and $\sgp_k(G_{r,s})\geq M$.

Finally, let $M=r+i(G_{r,s})$. This means that $r+i(G_{r,s})>k$, for otherwise we are in the earlier situation $M=k$. If $i(G_{r,s})=0$, then $M=r$, $A=Q$ and the deduction is trivial because vertices of $Q$ induce a clique. Suppose now that $A$ contains all $r$ vertices of $Q$ and first $i(G_{r,s})$ vertices of $I$ ordered by descending degree sequence. Let $B$ be any $k$-subset of $A$. First notice that $i(G_{r,s})<k$. Otherwise, if $i(G_{r,s})\geq k$, then $\delta(v_k)>r-k+k=r$, which is not possible as every vertex from $I$ can have at most $r$ neighbors. This means that $B_Q=B\cap Q\neq \emptyset$ and there are at least $k-i(G_{r,s})$ vertices from $Q$ in $B$. Since $\delta(v_i)\geq\delta(v_{i(G_{r,s})})>r-k+i(G_{r,s})$ for every $i\in[i(G_{r,s})]$, every vertex from $B_I=B\cap I$ contains more than $r-k+i(G_{r,s})$ neighbors in $Q$. In other words, every vertex from $B_I$ has a neighbor in $B_Q$. But then $B$ induces a connected subgraph and every Steiner $B$-tree contains only vertices from $B$. Therefore, $A$ is a $k$-Steiner general position set of $G_{r,s}$ and $\sgp_k(G_{r,s})\geq M$ follows.
\qed

\begin{figure}[ht!]
\centering
\begin{tikzpicture}[scale=.55, transform shape]
\node [draw, shape=circle,scale=0.7,fill=black] (a1) at  (0,0) {};
\node [draw, shape=circle,scale=0.7,fill=black] (a2) at  (0,1) {};
\node [draw, shape=circle,scale=0.7,fill=black] (a3) at  (0,2) {};
\node [draw, shape=circle,scale=0.7,fill=black] (a4) at  (0,3) {};
\node [draw, shape=circle,scale=0.7,fill=black] (a5) at  (0,4) {};
\node [draw, shape=circle,scale=0.7,fill=black] (a6) at  (0,5) {};
\node [draw, shape=circle,scale=0.7,fill=black] (a7) at  (0,6) {};
\node [draw, shape=circle,scale=0.7,fill=black] (a8) at  (0,7) {};
\node [draw, shape=circle,scale=0.7,fill=black] (a9) at  (0,8) {};

\node [draw, shape=circle,scale=0.7,fill=black] (b1) at  (1,0) {};
\node [draw, shape=circle,scale=0.7,fill=black] (b2) at  (1,1) {};
\node [draw, shape=circle,scale=0.7,fill=black] (b3) at  (1,2) {};
\node [draw, shape=circle,scale=0.7,fill=black] (b4) at  (1,3) {};
\node [draw, shape=rectangle,scale=1.2,fill=black] (b5) at  (1,4) {};
\node [draw, shape=circle,scale=0.7,fill=black] (b6) at  (1,5) {};
\node [draw, shape=circle,scale=0.7,fill=black] (b7) at  (1,6) {};
\node [draw, shape=circle,scale=0.7,fill=black] (b8) at  (1,7) {};
\node [draw, shape=circle,scale=0.7,fill=black] (b9) at  (1,8) {};

\node [draw, shape=circle,scale=0.7,fill=black] (c1) at  (2,0) {};
\node [draw, shape=circle,scale=0.7,fill=black] (c2) at  (2,1) {};
\node [draw, shape=circle,scale=0.7,fill=black] (c3) at  (2,2) {};
\node [draw, shape=rectangle,scale=1.2,fill=black] (c4) at  (2,3) {};
\node [draw, shape=circle,scale=0.7,fill=black] (c5) at  (2,4) {};
\node [draw, shape=circle,scale=0.7,fill=black] (c6) at  (2,5) {};
\node [draw, shape=circle,scale=0.7,fill=black] (c7) at  (2,6) {};
\node [draw, shape=circle,scale=0.7,fill=black] (c8) at  (2,7) {};
\node [draw, shape=circle,scale=0.7,fill=black] (c9) at  (2,8) {};

\node [draw, shape=circle,scale=0.7,fill=black] (d1) at  (3,0) {};
\node [draw, shape=circle,scale=0.7,fill=black] (d2) at  (3,1) {};
\node [draw, shape=rectangle,scale=1.2,fill=black] (d3) at  (3,2) {};
\node [draw, shape=circle,scale=0.7,fill=black] (d4) at  (3,3) {};
\node [draw, shape=circle,scale=0.7,fill=black] (d5) at  (3,4) {};
\node [draw, shape=circle,scale=0.7,fill=black] (d6) at  (3,5) {};
\node [draw, shape=circle,scale=0.7,fill=black] (d7) at  (3,6) {};
\node [draw, shape=circle,scale=0.7,fill=black] (d8) at  (3,7) {};
\node [draw, shape=circle,scale=0.7,fill=black] (d9) at  (3,8) {};

\node [draw, shape=circle,scale=0.7,fill=black] (e1) at  (4,0) {};
\node [draw, shape=rectangle,scale=1.2,fill=black] (e2) at  (4,1) {};
\node [draw, shape=circle,scale=0.7,fill=black] (e3) at  (4,2) {};
\node [draw, shape=circle,scale=0.7,fill=black] (e4) at  (4,3) {};
\node [draw, shape=circle,scale=0.7,fill=black] (e5) at  (4,4) {};
\node [draw, shape=circle,scale=0.7,fill=black] (e6) at  (4,5) {};
\node [draw, shape=circle,scale=0.7,fill=black] (e7) at  (4,6) {};
\node [draw, shape=rectangle,scale=1.2,fill=black] (e8) at  (4,7) {};
\node [draw, shape=circle,scale=0.7,fill=black] (e9) at  (4,8) {};

\node [draw, shape=circle,scale=0.7,fill=black] (f1) at  (5,0) {};
\node [draw, shape=circle,scale=0.7,fill=black] (f2) at  (5,1) {};
\node [draw, shape=circle,scale=0.7,fill=black] (f3) at  (5,2) {};
\node [draw, shape=circle,scale=0.7,fill=black] (f4) at  (5,3) {};
\node [draw, shape=circle,scale=0.7,fill=black] (f5) at  (5,4) {};
\node [draw, shape=circle,scale=0.7,fill=black] (f6) at  (5,5) {};
\node [draw, shape=rectangle,scale=1.2,fill=black] (f7) at  (5,6) {};
\node [draw, shape=circle,scale=0.7,fill=black] (f8) at  (5,7) {};
\node [draw, shape=circle,scale=0.7,fill=black] (f9) at  (5,8) {};

\node [draw, shape=circle,scale=0.7,fill=black] (g1) at  (6,0) {};
\node [draw, shape=circle,scale=0.7,fill=black] (g2) at  (6,1) {};
\node [draw, shape=circle,scale=0.7,fill=black] (g3) at  (6,2) {};
\node [draw, shape=circle,scale=0.7,fill=black] (g4) at  (6,3) {};
\node [draw, shape=circle,scale=0.7,fill=black] (g5) at  (6,4) {};
\node [draw, shape=rectangle,scale=1.2,fill=black] (g6) at  (6,5) {};
\node [draw, shape=circle,scale=0.7,fill=black] (g7) at  (6,6) {};
\node [draw, shape=circle,scale=0.7,fill=black] (g8) at  (6,7) {};
\node [draw, shape=circle,scale=0.7,fill=black] (g9) at  (6,8) {};

\node [draw, shape=circle,scale=0.7,fill=black] (h1) at  (7,0) {};
\node [draw, shape=circle,scale=0.7,fill=black] (h2) at  (7,1) {};
\node [draw, shape=circle,scale=0.7,fill=black] (h3) at  (7,2) {};
\node [draw, shape=circle,scale=0.7,fill=black] (h4) at  (7,3) {};
\node [draw, shape=rectangle,scale=1.2,fill=black] (h5) at  (7,4) {};
\node [draw, shape=circle,scale=0.7,fill=black] (h6) at  (7,5) {};
\node [draw, shape=circle,scale=0.7,fill=black] (h7) at  (7,6) {};
\node [draw, shape=circle,scale=0.7,fill=black] (h8) at  (7,7) {};
\node [draw, shape=circle,scale=0.7,fill=black] (h9) at  (7,8) {};

\node [draw, shape=circle,scale=0.7,fill=black] (i1) at  (8,0) {};
\node [draw, shape=circle,scale=0.7,fill=black] (i2) at  (8,1) {};
\node [draw, shape=circle,scale=0.7,fill=black] (i3) at  (8,2) {};
\node [draw, shape=circle,scale=0.7,fill=black] (i4) at  (8,3) {};
\node [draw, shape=circle,scale=0.7,fill=black] (i5) at  (8,4) {};
\node [draw, shape=circle,scale=0.7,fill=black] (i6) at  (8,5) {};
\node [draw, shape=circle,scale=0.7,fill=black] (i7) at  (8,6) {};
\node [draw, shape=circle,scale=0.7,fill=black] (i8) at  (8,7) {};
\node [draw, shape=circle,scale=0.7,fill=black] (i9) at  (8,8) {};

\draw(a1)--(a9);\draw(b1)--(b9);\draw(c1)--(c9);\draw(d1)--(d9);\draw(e1)--(e9);
\draw(f1)--(f9);\draw(g1)--(g9);\draw(h1)--(h9);\draw(i1)--(i9);
\draw(a1)--(i1);\draw(a2)--(i2);\draw(a3)--(i3);\draw(a4)--(i4);\draw(a5)--(i5);
\draw(a6)--(i6);\draw(a7)--(i7);\draw(a8)--(i8);\draw(a9)--(i9);

\draw[dotted](0,8)--(0,9);\draw[dotted](1,8)--(1,9);\draw[dotted](2,8)--(2,9);\draw[dotted](3,8)--(3,9);\draw[dotted](4,8)--(4,9);
\draw[dotted](5,8)--(5,9);\draw[dotted](6,8)--(6,9);\draw[dotted](7,8)--(7,9);\draw[dotted](8,8)--(8,9);
\draw[dotted](0,0)--(0,-1);\draw[dotted](1,0)--(1,-1);\draw[dotted](2,0)--(2,-1);\draw[dotted](3,0)--(3,-1);\draw[dotted](4,0)--(4,-1);
\draw[dotted](5,0)--(5,-1);\draw[dotted](6,0)--(6,-1);\draw[dotted](7,0)--(7,-1);\draw[dotted](8,0)--(8,-1);
\draw[dotted](-1,0)--(0,0);\draw[dotted](-1,1)--(0,1);\draw[dotted](-1,2)--(0,2);\draw[dotted](-1,3)--(0,3);\draw[dotted](-1,4)--(0,4);
\draw[dotted](-1,5)--(0,5);\draw[dotted](-1,6)--(0,6);\draw[dotted](-1,7)--(0,7);\draw[dotted](-1,8)--(0,8);
\draw[dotted](8,0)--(9,0);\draw[dotted](8,1)--(9,1);\draw[dotted](8,2)--(9,2);\draw[dotted](8,3)--(9,3);\draw[dotted](8,4)--(9,4);
\draw[dotted](8,5)--(9,5);\draw[dotted](8,6)--(9,6);\draw[dotted](8,7)--(9,7);\draw[dotted](8,8)--(9,8);

\node [scale=1.2] at (0,-1) {$-4$};\node [scale=1.2] at (1,-1) {$-3$};\node [scale=1.2] at (2,-1) {$-2$};
\node [scale=1.2] at (3,-1) {$-1$};\node [scale=1.2] at (4,-1) {$0$};\node [scale=1.2] at (5,-1) {$1$};
\node [scale=1.2] at (6,-1) {$2$};\node [scale=1.2] at (7,-1) {$3$};\node [scale=1.2] at (8,-1) {$4$};

\node [scale=1.2] at (-1.2,0) {$-4$};\node [scale=1.2] at (-1.2,1) {$-3$};\node [scale=1.2] at (-1.2,2) {$-2$};
\node [scale=1.2] at (-1.2,3) {$-1$};\node [scale=1.2] at (-1.2,4) {$0$};\node [scale=1.2] at (-1.2,5) {$1$};
\node [scale=1.2] at (-1.2,6) {$2$};\node [scale=1.2] at (-1.2,7) {$3$};\node [scale=1.2] at (-1.2,8) {$4$};

\end{tikzpicture}
\caption{By taking $k=4$, the vertices of the set $S$ appear squared.}
\label{fig:grid-k-4}
\end{figure}
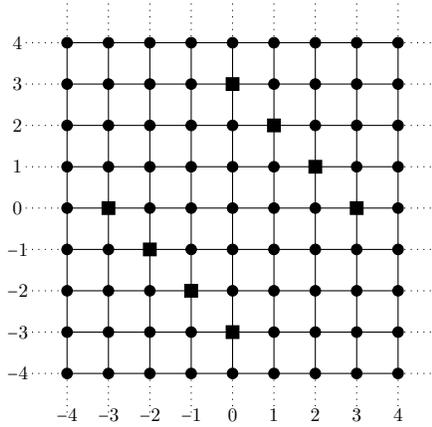

With $P_\infty$ we denote the two ways infinite path.
Let $V(P_{\infty})=\{\dots,-2,-1,0,1,2,\dots\}$ where $i$ is adjacent to $j$ if and only if $|i-j|=1$. The infinite grid $P_{\infty}\,\Box\, P_{\infty}$ is the Cartesian product of two infinite paths, that is $V(P_{\infty}\,\Box\, P_{\infty})=\{(i,j)\,:\,i,j\in \mathbb{Z}\}$ and $(i,j)(k,\ell)\in E(P_{\infty}\,\Box\, P_{\infty})$ when $|i-j|+|k-\ell|=1$, see Fig.~\ref{fig:grid-k-4}.

\begin{theorem}\label{grid}
$\sgp_k(P_{\infty}\,\Box\, P_{\infty})\ge 2k$.
\end{theorem}

\proof
Let $S=S_1\cup S_2$ be a set of vertices with $S_1= \{(k-1,0),(k-2,1),\dots,(1,k-2),(0,k-1)\}$ and $S_2=\{(-k+1,0),(-k+2,-1),\dots,(-1,-k+2),(0,-k+1)\}$. See Fig.~\ref{fig:grid-k-4} for an example when $k=4$. Notice that $|S_1|=|S_2|=k$. We will show prove that $S$ is a $k$-Steiner general position set of $P_{\infty}\,\Box\, P_{\infty}$. Let $B\subset S$ be a $k$-subset. If $B=S_1$ or $B=S_2$, then it can be easily noted that any Steiner $B$-tree does not contain any vertex of $S\setminus B$. Hence, we may assume $B_1=B\cap S_1\ne \emptyset$ and $B_2=B\cap S_2\ne \emptyset$.

\begin{figure}[ht!]
\centering
\begin{tikzpicture}[scale=.6, transform shape]

\draw(0,-4)--(0,4);
\draw(-4.5,0)--(4,0);

\node [draw, shape=circle,scale=0.7,fill=black] (a) at  (-4,-2) {};
\node [draw, shape=circle,scale=0.7,fill=black] (b) at  (-2.5,-3.5) {};
\node [draw, shape=circle,scale=0.7,fill=black] (c) at  (-2.5,-2) {};

\node [draw, shape=circle,scale=0.7,fill=black] (d) at  (3,2) {};
\node [draw, shape=circle,scale=0.7,fill=black] (e) at  (1.5,3.5) {};
\node [draw, shape=circle,scale=0.7,fill=black] (f) at  (1.5,2) {};

\draw[thick] (c) rectangle (f);
\draw[thick] (a)--(c)--(b);
\draw[thick] (d)--(f)--(e);
\draw[dotted] (d)--(3,0);
\draw[dotted] (a)--(-4,0);
\draw[dotted] (b)--(0,-3.5);
\draw[dotted] (e)--(0,3.5);

\node [scale=1.2] at (0,-4.4) {$0$};
\node [scale=1.2] at (-5,0) {$0$};
\node [scale=1.2] at (3.2,-0.3) {$i_2$};
\node [scale=1.2] at (1.85,-0.3) {$i_1$};
\node [scale=1.2] at (-4,0.3) {$i_4$};
\node [scale=1.2] at (-2.15,0.3) {$i_3$};
\node [scale=1.2] at (0.35,-3.5) {$j_3$};
\node [scale=1.2] at (0.35,-1.7) {$j_4$};
\node [scale=1.2] at (-0.35,3.5) {$j_1$};
\node [scale=1.2] at (-0.35,1.7) {$j_2$};

\end{tikzpicture}
\begin{tikzpicture}[scale=.6, transform shape]

\draw(0,-4)--(0,4);
\draw(-4.5,0)--(4,0);

\node [draw, shape=circle,scale=0.7,fill=black] (a) at  (-4,-2) {};
\node [draw, shape=circle,scale=0.7,fill=black] (b) at  (-2.5,-3.5) {};
\node [draw, shape=circle,scale=0.7,fill=black] (c) at  (-2.5,-2) {};

\node [draw, shape=circle,scale=0.7,fill=black] (d) at  (1.5,2) {};
\node [draw, shape=circle,scale=0.7,fill=black] (e) at  (1.5,2) {};
\node [draw, shape=circle,scale=0.7,fill=black] (f) at  (1.5,2) {};

\draw[thick] (c) rectangle (f);
\draw[thick] (a)--(c)--(b);
\draw[thick] (d)--(f)--(e);
\draw[dotted] (a)--(-4,0);
\draw[dotted] (b)--(0,-3.5);

\node [scale=1.2] at (0,-4.4) {$0$};
\node [scale=1.2] at (-5,0) {$0$};
\node [scale=1.2] at (2.35,-0.3) {$i_1=i_2$};
\node [scale=1.2] at (-4,0.3) {$i_4$};
\node [scale=1.2] at (-2.15,0.3) {$i_3$};
\node [scale=1.2] at (0.35,-3.5) {$j_3$};
\node [scale=1.2] at (0.35,-1.7) {$j_4$};
\node [scale=1.2] at (-0.85,2.3) {$j_2=j_1$};

\end{tikzpicture}
\caption{Two sketches of possibilities for the set $Y$. In the second one, $i_1=i_2$ and $j_2=j_1$, or equivalently, $(i_1,j_1)=(i_2,j_2)$.}\label{fig:Steiner-tree-shape}
\end{figure}
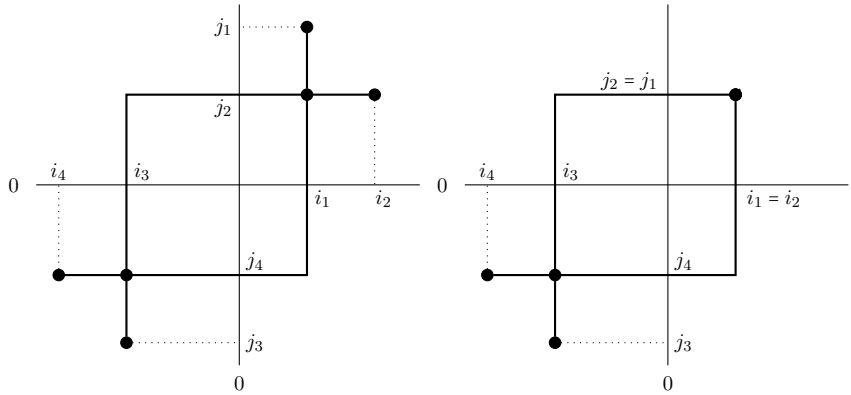

Let $j_1,i_2$ and $j_2,i_1$ be the largest and the smallest indexes, respectively, such that the vertices $(i_1,j_1),(i_2,j_2)$ belong to $S_1$. Analogously, let $j_4,i_3$ and $j_3,i_4$ be the largest and the smallest indexes, respectively, such that the vertices $(i_3,j_3),(i_4,j_4)$ belong to $S_2$. We remark that it could happen $(i_1,j_1)=(i_2,j_2)$ or $(i_3,j_3)=(i_4,j_4)$. Let $X=\{(i_1,j_1),(i_2,j_2),(i_3,j_3),(i_4,j_4)\}$ and consider a Steiner $X$-tree. According to the structure of the set $S$, we notice that only those vertices of $Y$ given below could belong to a Steiner $X$-tree. See Fig.~\ref{fig:Steiner-tree-shape} for a sketch of $Y$.
\begin{align*}
  Y= & \left(\{i_3,\dots,i_1\}\times \{j_4,\dots,j_2\}\right)\cup \left(\{i_1,\dots,i_2\}\times\{j_2\}\right)\cup \left(\{i_1\}\times\{j_2,\dots,j_1\}\right) \\
   &\cup \left(\{i_4,\dots,i_3\}\times\{j_4\}\right)
\cup \left(\{i_3\}\times\{j_3,\dots,j_4\}\right).
\end{align*}

Now, we first observe that every vertex $(i,j)\in S$ such that $i\notin [i_1:i_2]\cup [i_4:i_3]$ and $j\notin [j_2:j_1]\cup [j_3:j_4]$ does not belong to $Y$, and so, it is not included in any Steiner $X$-tree.

On the other hand, if $(i,j)\in B_1\setminus X$, then to obtain a Steiner $X\cup \{(i,j)\}$-tree we can only add to an Steiner $X$-tree some vertices of the set $\{i_1,\dots,i\}\times \{j_2,\dots,j\}$. This allows to claim that if $(i',j')\notin B_1$, then it will not belong to any Steiner $X\cup \{(i,j)\}$-tree for every $(i,j)\in B_1\setminus X$. This procedure can be iterated for all the remaining vertices of $B_1$, and consequently, we will obtain that any Steiner $X\cup B_1$-tree does not contain vertices from $S\setminus B_1$. By using some symmetrical arguments, we will obtain that any Steiner $X\cup (B_1\cup B_2)$-tree does not contain vertices from $S\setminus (B_1\cup B_2)$, which is precisely that any Steiner $B$-tree does not contain vertices from $S\setminus B$. Therefore, $S$ is a $k$-Steiner general position set for $P_{\infty}\,\Box\, P_{\infty}$, and the lower bound follows.
\qed

\section{Concluding remarks and problems}

In Proposition~\ref{prop:k-1-and-not-k} we have demonstrated that a $k$-Steiner general position need not be  a $k'$-Steiner general position set for $k' > k$. This result indicates that there is no monotony relation for $k$-Steiner general position sets with respect to inclusion, for every graph $G$, but not necessarily for the value of the parameter $\sgp_k(G)$. Hence, we pose the following problem.

\begin{problem}
Is there any monotony relation between $\sgp_k(G)$ and $\sgp_{k+1}(G)$?
\end{problem}

It is already known that computing the general position number of graphs is NP-hard in general. In this sense, the answer to the following problem seems obvious.

\begin{problem}
Determine the complexity of computing the $k$-Steiner general position number.
\end{problem}

We are not aware of any lexicographic product for which the bound of Theorem~\ref{lex} is not sharp, hence se pose:

\begin{problem}
Is the bound of Theorem~\ref{lex} sharp for all lexicographic products?
\end{problem}

It is easy to construct several split graphs such that the bound of Theorem~\ref{split} is sharp. It remains to describe all split graphs for which the equality holds.

\begin{problem}
For which split graphs is the bound of Theorem~\ref{split} sharp?
\end{problem}

The bound from Theorem~\ref{grid} is tight because $\sgp_2(P_{\infty}\,\Box\, P_{\infty})=\gp(P_{\infty}\,\Box\, P_{\infty})=4$ was proved in~\cite{manuel-2018b}. We wonder whether a parallel result holds for each $k > 2$: 

\begin{problem}
Does the equality $\sgp_k(P_{\infty}\,\Box\, P_{\infty})=2k$ holds for $k>2$?
\end{problem}

Finally, in~\cite{klavzar-2021} the general position number of arbitrary integer lattices was determined, that is, of the Cartesian product of finitely many factors $P_{\infty}$. Hence we also pose: 

\begin{problem}
Investigate $\sgp_k(P_{\infty}\,\Box\, \cdots \,\Box\, P_{\infty})$ for $k>2$. 
\end{problem}

\section*{Acknowledgements}

Sandi Klav\v{z}ar and Iztok Peterin acknowledge the financial support from the Slovenian Research Agency (research core funding No.\ P1-0297 and projects J1-9109, J1-1693, N1-0095, N1-0108). Dorota Kuziak and Ismael G. Yero have been partially supported by the Spanish Ministry of Science and Innovation through the grant PID2019-105824GB-I00.

\end{document}